                       \documentstyle[12pt,amssymb,psfig]{amsart}
                        \numberwithin{equation}{section}
                        \theoremstyle{plain}
                         \newtheorem{thm}{Theorem}[section]
                                \newtheorem{lem}[thm]{Lemma}
                                \newtheorem{pro}[thm]{Proposition}
                                \newtheorem{cor}[thm]{Corollary}
                        \theoremstyle{definition}
                                \newtheorem{rem}{Remark}[thm]

\newcommand{\psdraw}[2]
         {\begin{array}{c} \hspace{-1.3mm}
         \raisebox{-4pt}{\psfig{figure=#1.eps,width=#2}}
         \hspace{-1.9mm}\end{array}}

\newcommand{\hh}{{\frak h}}
\newcommand{\al}{\alpha}

\newcommand{\tr}{\text{\rm tr}}

\newcommand{\Q}{{\Bbb Q}}
\newcommand{\C}{{\Bbb C}}

\newcommand{\Wr}{\Weyl_r}
\newcommand{\tWr}{\hat{\Weyl}_r}

\newcommand{\R}{{\Bbb R}}

\newcommand{\Z}{{\Bbb Z}}

\newcommand{\fg}{{\frak g}}

\newcommand{\UU}{{\cal U}}

\newcommand{\cCfin}{{\cal C}}

\newcommand{\Weyl}{W}
\newcommand{\sn}{{\rm sn}}
\newcommand{\id}{{\rm id}}

\newcommand{\bCr}{\bar C_r}

\begin{document}

                \title[Quantum invariants of 3-manifolds]
{Quantum invariants of 3-manifolds: integrality, splitting, and perturbative
                expansion}

                \author[ Thang Le]{Thang T. Q. Le }
\thanks{This work is partially supported by NSF grant DMS-9626404.}
\address{Dept. of Mathematics, SUNY at Buffalo, Buffalo, NY 14214, USA }
\email{letu@@math.buffalo.edu}
                \maketitle
\begin{abstract}
We consider quantum invariants of 3-manifolds associated with arbitrary simple Lie
algebras. Using the symmetry principle we show how to decompose the quantum invariant as
the product of two invariants, one of them is the invariant corresponding to the
projective group. We then show that the projective quantum invariant is always an
algebraic integer, if the quantum parameter is a prime root of unity. We also show that
the projective quantum invariant of rational homology 3-spheres has a perturbative
expansion a la Ohtsuki. The presentation of the theory of quantum 3-manifold is 
self-contained. 
\end{abstract}

\addtocounter{section}{-1}

\section{Introduction}

\subsection{} For a simple Lie algebra $\fg$ over $\C$ with Cartan matrix $(a_{ij})$ let
$d = \max_{i\neq j}|a_{ij}|$. Thus $d=1$ for the $ADE$ series, $d=2$ for $BCF$ and $d=3$
for $G_2$. The quantum group associated with $\fg$ is a Hopf algebra over $\Q(q^{1/2})$,
where $q^{1/2}$ is the quantum parameter. To fix the order let us point out that our $q$
is $q^2$ in \cite{Kassel,Kirillov,Turaev} or $v^2$ in the book \cite{Lusztig}. For
example, the quantum integer is given by

$$[n] =\frac{q^{n/2}-q^{-n/2}} {q^{1/2}-q^{-1/2}}.$$

\subsection{} Modular categories, and hence quantum 3-manifold invariants associated with
$\fg$, can be defined only when $q$ a root of unity of order $r$ divisible by $d$, since
this fact guarantees that the so-called S-matrix is invertible. For $r$ not divisible by
$d$, quantum invariants can still be defined, but modular categories might not exist. In
this paper we will focus mainly in the more general situation, when $r$ may or may not be
divisible by $d$. The reason is eventually we want $r$ to be a prime number. Note that
3-manifold invariants for the case when $r$ is not divisible by $d\neq 1$ have not been
studied earlier.

We will present a self-contained theory of quantum 3-manifold invariants, for arbitrary 
simple Lie algebra. By making use of an integrality result (Proposition \ref{s10}) we 
will establish the existence of quantum invariants without using the theory of quantum 
groups at roots of unity. (We use quantum group and link invariants with general 
parameter, and only in the last minute, replace $q$ by a root of unity.) The quantum 
invariant of a manifold $M$ will be denoted by $\tau^\fg_M(q)$, considered as a function 
with domain roots of unity. 

\subsection{} Although the usual construction of modular category might fail, say when $r$ is
a prime number and $d\neq 1$, we will show  that by using the root lattice instead of the
weight lattice in the construction, one can still get a modular category. Actually the
construction goes through for a much larger class of numbers $r$ -- one needs only that
$r$ is coprime with $d \det (a_{ij})$. The corresponding 3-manifolds invariant, denoted
by $\tau_M^{P\fg}(q)$, could be considered as the invariant associated to the projective
group -- the smallest complex Lie group whose Lie algebra is $\fg$. The reason is that
the set of all highest weights of modules of the projective group spans the root lattice.
As in the case of $\fg$, the invariant $\tau^{P\fg}_M$ can also  be defined  when $r$ is
not coprime with $d \det (a_{ij})$ (although modular categories might not exist).

\subsection{} We will show that in most cases, $\tau_M^{P\fg}$ is finer than
$\tau^\fg_M$. More precisely, if $r$ is coprime with $\det(a_{ij})$, then

\begin{equation} \label{first}
\tau_M^\fg = \tau_M^{P\fg} \times \tau^G_M,
\end{equation}
where $\tau^G_M$ is the 3-manifold invariant associated with the center group $G$ (which
is isomorphic to the quotient of the weight lattice by the root lattice) and a naturally
defined bilinear form on it. The invariant $\tau^G_M$ is a weak invariant, since it is 
determined by the first homology group and the linking form on its torsion (see 
\cite{MOO,Deloup,Turaev}). In some cases $\tau^G_M\equiv 0$, and hence $\tau^\fg$ is 
trivial, although $\tau^{P\fg}$ is not. 

The invariant associated to the projective group was first introduced and the splitting
(\ref{first}) was obtained by Kirby and Melvin for $\fg=sl_2$, and Kohno and Takata for
$\fg=sl_n$. Recently Sawin \cite{Sawin}, based on the work of M\"uger and Bruguierre,
established a similar result, but he considered only the case $r$ divisible by $d$, and
the group $G$ cyclic (so he excludes a half of the series $D$ case). When $d=2$, the
result (\ref{first}) complements Sawin's work, i.e. it covers the case that is not
considered by Sawin. When $d=1$, (\ref{first}) overlaps with Sawin's work. But even in 
this case, our method is quite different, it can be uniformly applied to any simple Lie 
algebra, and in addition, we get the integrality and perturbative expansion of quantum 
3-manifold invariants (see below).

\subsection{} We will show that unlike $\tau^\fg$, the ``projective" invariant
$\tau^{P\fg}_M$ is always an algebraic integer, provided that the order $r$ is an odd
prime. In fact, we will prove that in this case, $\tau^{P\fg}_M(q) \in \Z[q]$. Apriori,
both $\tau^{P\fg}_M$ and $\tau^\fg_M$ are rational functions in a fractional power of
$q$. Integrality of $\tau^{P\fg}_M$ for $\fg=sl_2$ was first established by H. Murakami
\cite{Murakami} by difficult computations, for $\fg =sl_n$ by Takata-Yokota \cite{TY} and
Masbaum-Wenzl \cite{MW}, based on an idea of Roberts.  We will use a different approach 
that is good for all simple Lie algebras. 

\subsection{}  Finally we will show that $\tau^{P\fg}_M$, with $M$ a rational homology
3-sphere, has a ``perturbative expansion a la Ohtsuki". The function $\tau^{P\fg}_M$ can
be defined only at roots of unity, and we want to expand it around $q=1$. For the case
$\fg =sl_2$, Ohtsuki showed that there exists a kind of number-theoretic expansion, which
we call perturbative expansion. We proved a similar result for $\fg = sl_n$ in
\cite{perturbative} and will extend the result to other Lie algebras here. We borrowed an
idea using Gauss integrals from Rozansky's work \cite{Rozansky}, although we will not
explicitly use Gauss integral.

\subsection{}
 The paper is organized as
follows. In section 1 we recall  quantum link invariants and their important properties:
integrality and symmetry at roots of 1.  In section 2 we present the general theory of
3-manifold quantum invariants (not using the theory of quantum groups at roots of unity).
Invariants associated to the projective group are considered in section 3. Their
integrality is proved in section 4. Section 5 is devoted to the perturbative expansion.

The paper is part of the talks the author gave at conference ``Knots in Hellas" (Delphi, 
Greece, July 1998),  Borel seminars (Bern university, June 1999), Summer School on 
quantum invariants (Grenoble, France, July 1999), and Conference on quantum invariants 
(Calgary, July 1999). The author would like to thank the organizers, especially J. 
Przytycki,  N. Habegger, C. Lescop, and J. Bryden for inviting him to give talks at these 
conferences. He would like to thank G. Masbaum, H. Murakami and V. Turaev for helpful 
discussions, and the Mittag-Leffler Institute for hospitality and support during May 
1999.

\section{Quantum link invariants: Integrality and Symmetries}

\subsection{Lie algebras and Quantum groups}
We recall here some facts from the theory of Lie algebras and quantum groups, mainly in
order to fix notation. For the theory of quantum groups, see \cite{Kassel,Lusztig}. 

\subsubsection{Lie algebra}
Let $(a_{ij})_{1\le i,j\le \ell}$ be the Cartan matrix of a simple complex Lie algebra
$\fg$. There are relatively prime integers $d_1,\dots,d_\ell$ in $\{1,2,3\}$ such that
the matrix $(d_ia_{ij})$ is symmetric. Let $d$ be the maximal of $(d_i)$. The values of
$d$, and other data, for various Lie algebras are listed in Table 1.

\begin{table}[ht]
\begin{center}
\begin{tabular}{|c|c|c|c|c|c|c|c|c|c|c|c|}
\hline
     & $A_\ell$ & $B_\ell$    &  $B_\ell$  & $C_\ell$   & $D_\ell$  &  $D_\ell$  & $E_6$  &  $E_7$  &  $E_8$  & $F_4$  & $G_2$ \\
     &          & $\ell$ odd  & $\ell$ even &            & $\ell$ odd&$\ell$ even &        &         &         &        &       \\ \hline
$d$  & $1$      & $2$         & 2           &  2         &  1        &   1        &   1    &    1    &    1    &   2    &   3   \\ \hline
$D$  & $\ell+1$ & 2           &  1          &    1       &   4       &  2         &   3    &   2     &   1     &   1    &   1   \\ \hline
$G$  &$\Z_{\ell+1}$& $\Z_2$   & $\Z_2$      & $\Z_2$     & $\Z_4$    &$\Z_2\times\Z_2$&$\Z_3$&$\Z_2$ &   1     &  1     & 1     \\ \hline
$h$  & $\ell+1$ & $2\ell$     & $2\ell$     & $2\ell$    &$2\ell-2$  &$2\ell-2$   & 12     &  18     & 30      &   12   & 6     \\ \hline
$ h^\vee$
     & $\ell +1$& $2\ell-1$   & $2\ell-1$   & $\ell+1$   &$2\ell-2$  &$2\ell-2$   & 12     &  18     & 30      &   9    &   4  \\ \hline
\end{tabular}
\end{center}
\caption{}
\end{table}

We fix a Cartan subalgebra $\hh$ of $\fg$ and  basis roots $\al_1,\dots,\al_\ell$ in the
dual space $\hh^*$. Let $\hh^*_\R$ be the $\R$-vector space spanned by
$\al_1,\dots,\al_\ell$. The root lattice $Y$ is the $\Z$-lattice  generated by
$\al_i,i=1,\dots,\ell$. Define the scalar product  on $\hh^*_\R$ so that $(\al_i| \al_j)=
d_i a_{ij}$.
 Then  $(\alpha|\alpha)=2$ for every {\em short} root $\alpha$.

Let $\Z_+$ be the set of all non-negative integers. The weight lattice  $X$ (resp. the
set of dominant weights $X_+$) is the set of all $\lambda \in\hh^*_\R$ such that $\langle
\lambda, \al_i \rangle := \frac{2(\lambda|\al_i)}{(\al_i|\al_i)} \in\Z$ (resp.  $\langle
\lambda, \al_i \rangle\in\Z_+$) for $i=1,\dots,\ell$. Let $\lambda_1,\dots,\lambda_\ell$
be the fundamental weights, i.e. the $\lambda_i\in \hh^*_\R$ are defined by $\langle
\lambda_i,\al_j\rangle =\delta_{ij}$, or $(\lambda_i|\al_j)=d_i \delta_{ij}$. Then $X$ is
the $\Z$-lattice generated by $\lambda_1,\dots, \lambda_\ell$. The root lattice $Y$ is a
subgroup of the weight lattice $X$, and the quotient $G=X/Y$ is called the {\em
fundamental group}. If $\mu\in X$ and $\al\in Y$, then $(\mu|\al)$ is always an integer.
On the root lattice $Y$, the form $(\cdot|\cdot)$ is even.

Let $\rho$ be  the half-sum of all positive roots. Then $\rho=\lambda_1+\dots
+\lambda_\ell \in X_+$, and $2\rho \in Y$.

Let $C$ denote the fundamental chamber:

$$ C= \{ x\in \hh^*_\R \mid (x|\al_i)\ge 0, \quad i=1,\dots,\ell\}.$$

The Weyl group $\Weyl$ is the group generated by reflections in the walls of $C$. In the
chamber $C$ there is exactly one root of length $\sqrt2$; it is called the short highest
root, and denoted by $\alpha_0$.

For a positive integer $r$ let

$$C_r = \{ x\in C \mid (x|\al_0) < r\}.$$

Then the topological closure $\bar C_r$ is a simplex. The reflections in the walls of
$\bar C_r$ generate the affine Weyl group $\Weyl_r$. One also has $ \Weyl_r = W \ltimes
rY,$ where $rY$ denotes the translations by vectors $ry, y\in Y$.

Finite-dimensional simple $\fg$-modules are parametrized by $X_+$: for every
$\lambda\in X_+$, there corresponds a unique simple $\fg$-module $\bar \Lambda_\lambda$.

The Coxeter number and the dual Coxeter numbers are defined by $ h = 1 + (\al_0|\rho)$
and $h^\vee = 1 + \max_{\al > 0} \frac{(\al|\rho)}{d}$, see Table 1. Note that $d h^\vee
\ge h$.

\subsubsection{The quantum group  $\UU$ and its category of representations}
The quantum group $\UU=\UU_q(\fg)$ associated to $\fg$ is a Hopf algebra defined over
$\Q[q^{\pm 1/2D}]$, where $q^{1/2}$ is the quantum parameter and $D$ is the least 
positive integer such that $(\mu|\mu')\in \frac{1}{D}\Z$ for every $\mu,\mu'\in X$ (see 
\cite{Lusztig}). The category $\cCfin$ of finite-dimensional $\UU$-modules of type 1 is a 
ribbon category. In this paper we consider only $\UU$-modules of type 1. The introduction 
of the fractional power $q^{1/2D}$ is necessary for the definition of the braiding. 
Finite-dimensional simple $\UU$-modules of type 1 are also parametrized by $X_+$: for 
every $\lambda\in X_+$, there corresponds a unique simple $\UU$-module $ 
\Lambda_\lambda$, a deformation of $ \bar \Lambda_\lambda$. Actually, the Grothendieck 
ring of finite-dimensional $\UU$-modules of type 1 is isomorphic to that of 
finite-dimensional $\fg$-modules. 

The reader should not confuse our $q$ with the quantum parameter used in the definition
of quantum groups by several authors. For example, our $q$ is equal to $q^2$ in
\cite{Kassel,Kirillov,Turaev}, or $v^2$ in Lusztig book \cite{Lusztig}.

\subsection{Quantum link invariants}

\subsubsection{General}

Suppose $L$ is a framed oriented link with $m$ ordered components, then the quantum
invariant $J_L(V_1,\dots,V_m)$, for $V_1,\dots,V_m\in \cCfin$, is defined (since $\cal C$
is a ribbon category), with values in $\Z[q^{1/2D}]$.  The modules $V_1,\dots,V_m$ are
usually called the colors. The fact that $J_L$ has integer coefficients follows from
Lusztig's theory of canonical basis (see a detailed proof in \cite{integrality}).

We will also use another normalization of the quantum invariant:

$$Q_L(V_1,\dots,V_m):= J_L(V_1,\dots,V_m) \times J_{U^{(m)}}(V_1,\dots,V_m).$$

Here $U^{(m)}$ is the 0 framing trivial link of $m$ components. This normalization is
more suitable for the study of quantum 3-manifold invariants, and will help us to get rid
of the $\pm$ sign in many formulas.

Since finite-dimensional irreducible $\UU$-modules are parametrized by $X_+$, we define

$$Q_L(\mu_1,\dots,\mu_m) := Q_L(\Lambda_{\mu_1-\rho},\dots,\Lambda_{\mu_m-\rho}).$$

Note the  shift by $\rho$. This definition is good only for $\mu_j\in \rho+ X_+= X \cap
(\text{interior of $C$})$. We define $Q_L(\mu_1,\dots,\mu_m)$ for arbitrary $\mu_j\in X$
by requiring that $Q_L(\mu_1,\dots,\mu_m)=0$ if one of the $\mu_j$'s is on the boundary
of $C$, and that $Q_L(\mu_1,\dots,\mu_m)$ is component-wise invariant under the action of
the Weyl group $W$, i.e. for every $w_1,\dots,w_m\in W$,

$$ Q_L(w_1(\mu_1),\dots,w_m(\mu_m)) = Q_L(\mu_1,\dots,\mu_m).$$

\subsubsection{Example} 
\begin{figure}[htpb]
$$ \psdraw{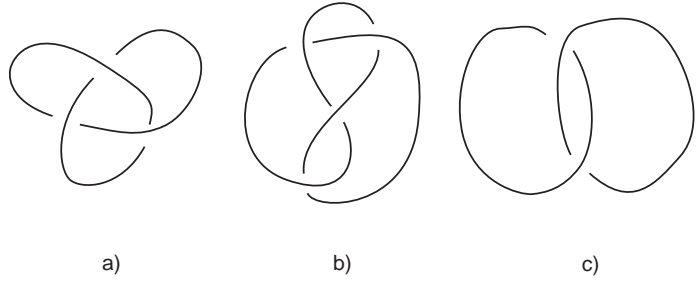}{3.6in} $$
\caption{The trefoil, figure 8 knots and the Hopf link}\label{trefoil}
\end{figure}

Suppose $\fg=sl_2$.  For a knot $K$, the invariant $J_K(N)$, with $N$ a positive
integer, is known as the colored Jones polynomial. Here $N$ stands for the unique simple
$sl_2$-module of dimension $N$. Suppose $K$ is the right-hand trefoil, see Figure
\ref{trefoil}. Then

$$ J_K(N) = [N]\, q^{1-N} \sum_{n=0}^\infty q^{-nN} (1-q^{1-N}) (1-q^{2-N}) \dots
(1-q^{n-N}).$$

The sum is actually finite, for any positive integer $N$.  Similar formulas have also
been obtained by Gelca and Habiro.

For the figure 8 knot (also obtained by Habiro)

$$  J_K(N) = [N] \sum_{n=0} ^\infty q^{-n N} (1-q^{N -1})(1-q^{N -2})
\dots (1-q^{N -n}) \times (1-q^{N +1})(1-q^{N +2}) \dots (1-q^{N+m}).$$ 

\subsubsection{The trivial knot}\label{trivial}
Suppose $U$ is the trivial knot. Then $J_U(V)$ is called the {\em quantum dimension} of
$V$; its value is well-known:
\begin{align}
J_U(\mu) &= \frac{1}{\psi}\prod_{\text {positive roots } \al}( q^{(\mu  | \al)/2} -
q^{-(\mu | \al)/2}) \label{unknot1}
\\
&=\frac{1}{\psi}\,{\sum _{w\in \Weyl} \sn(w) q^{ (\mu | w(\rho))}}
, \label{unknot2}
\end{align}
where  $\sn(w)$ is the sign of $w$ and
\begin{equation} \psi =  \prod_{\al >0} (q^{(\rho | \al)/2} - q^{-(\rho | \al)/2}) = \sum_{w\in W}
\sn(w) q^{(\rho|w(\rho))}.
\label{psi}
\end{equation}

Note that $J_U(\rho) =1$, and $J_U(-\mu) = (-1)^s J_U(\mu)$, where $s$ is the number of
positive roots. Also

\begin{equation}\psi = q^{-|\rho|^2}\prod_{\al >0}( q^{(\al|\rho)} -1). \label{ssss}
\end{equation}
 Hence if $q$ is a root of
unity of order $r \ge d h^\vee \ge 1 + \max_{\al >0} (\al|d)$, then $\psi \neq 0$.

\subsubsection{The Hopf Link}

Let $H$ be the Hopf link, see Figure \ref{trefoil},  with framing 0 on each component.
Then

\begin{equation}
 J_H(\mu,\lambda) = \frac{1}{\psi}\displaystyle{\sum _{w\in \Weyl} \sn(w) q^{ (\mu |
w(\lambda))}}. \label{Hopflink}
\end{equation}

For a proof, see \cite{Turaev}.  Note that $J_H(\mu,\rho) = J_U(\mu)$.

\subsection{Integrality}

In general, $J_L$ and $Q_L$  contain fractional powers of $q$. The integrality,
formulated below, shows that the fractional powers can be factored out.

\begin{thm} [Integrality, \cite{integrality}] Suppose $\mu_1,\dots,\mu_m \in X_+$. Then
  $Q_L(\Lambda_{\mu_1},\dots,\Lambda_{\mu_m})$ is in $ q^{p/2}
\Z[q^{\pm 1}]$,
where $p$ is a (generally fractional) number determined by the linking matrix $l_{ij}$ of
$L$:
$$p= \sum_{1\le i,j\le m} l_{ij} (\mu_i| \mu_j) + \sum_{1\le i\le m} l_{ii} (2\rho|
\mu_i) \quad \in \quad \frac{1}{D}\Z.$$
\end{thm}

If all $\mu_j$'s are in the {\em root lattice}, then the number $p$ is even. Hence we
have

\begin{cor}
If all the $\mu_j$'s are in the {\em root lattice}, then
$Q_L(\Lambda_{\mu_1},\dots,\Lambda_{\mu_m})$ is in $\Z[q^{\pm1}]$.
\label{integ}
\end{cor}

\subsection{The First Symmetry Principle}

\newcommand{\eJ}{{_\varepsilon J}}
\newcommand{\eqr}{\;\overset{{(r)}}{=}\;}

Suppose $f,g\in  \Z[q^{\pm 1/2D}]$. We say that $f$ equals $g$ {\em at  $r$-th roots of
unity} and write $$ f\eqr g,$$ if there is a number $a\in \frac{1}{2D}\Z$ such that $f,g
\in q^{a}\Z[q^{\pm 1}]$, and for every $r$-th root of unity $\xi$, one has $$q^{-a}
f|_{q=\xi} \quad  =\quad q^{-a} g|_{q=\xi}.$$ There is no need to fix a $2D$-th root of
$\xi$.

An equivalent definition: $f,g\in \Z[q^{\pm 1/2D}]$ are equal at $r$-th roots of unity if
for every  $2Dr$-th root $\zeta$ of unity, one has $f|_{q^{1/2D}=\zeta} \quad =\quad
g|_{q^{1/2D}=\zeta}$.

Recall that the simplex  $\bar C_r$ is a fundamental domain of the affine Weyl group
$\Wr$.

\begin{thm}[First Symmetry Principle, see \cite{integrality}]
At  $r$-th roots of unity, the quantum invariant $Q_L$ is {\em component-wise} invariant
under the action of the affine Weyl group $\Wr$. This means, for every $w_1,\dots,w_m\in
\Wr$,
\begin{equation}
 Q_L\left(w_1(\mu_1),\dots,w_m(\mu_m)\right) \eqr Q_L(\mu_1,\dots,\mu_m).
\notag
\end{equation}

If one of the $\mu_1,\dots,\mu_m$ is on the boundary of $\bCr$, then
$J_L(\mu_1,\dots,\mu_m)\eqr 0$.
\end{thm}

\subsection{The Second Symmetry Principle}
\subsubsection{Action of $G$ on $\bar C_r$}
There is an action of the center group $G= X/Y$ on $\bar C_r$, and although $Q_L$ is not
really invariant under this action, it almost is. Let us first describe the action.
 For $\mu\in \bar C_r$ and $g\in G = X/Y$ let
$\tilde g\in X$ be a lift of $g$, and define

$$ g(\mu) = \mu+r \tilde g   \in  (\text{$X$ \hspace {-.4cm} $\mod{W_r}$})\quad  \equiv \bar C_r.$$

Another way to look at the action of $G$ is the following. Recall that $\Wr= \Weyl\ltimes
rY$. Note that $X$ is  invariant under the action of the Weyl group. Let $\tWr$ be the
group generated by $\Weyl$ and translations by $rX$. Then $\tWr = \Weyl \ltimes rX$. If
$\lambda\in X$ and $w\in \Weyl$, then $w(\lambda)-\lambda$ is in $Y$. This implies  $\Wr$
is a {\em normal subgroup} of $\tWr$. We have an exact sequence
\begin{equation*}
 1\to \Wr \to  \tWr \to G \to 1.
\end{equation*}

Taking the action of $\tWr$ modulo the action of $\Wr$, we get the action of $G$ on $\bar
C_r$. For more details and examples of actions of $G$, see \cite{integrality}.

On the group $G$ there is a  symmetric bilinear form with values in $\Q/\Z$ defined as
follows. Suppose $g_1,g_2\in G=X/Y$. Let $\tilde g_1,\tilde g_2\in X$ be respectively
lifts of $g_1,g_2$. Define $(g_1|g_2) := (\tilde g_1|\tilde g_2) \in \Q/\Z$. (Actually,
$(g_1|g_2) \in \frac{1}{D}\Z/\Z$.)  Similarly, for $\mu\in X$ and $g\in G$ one can define
$(\mu|g)\in \frac{1}{D}\Z/\Z$.

\subsubsection{Second Symmetry Principle}

\begin{thm}(\cite{integrality})
\label{super} Suppose $\mu_1,\dots,\mu_m\in \bar C_r$ and
$g_1,\dots,g_m\in G$. Then
\begin{equation}
 Q_L\left(g_1(\mu_1),\dots,g_m(\mu_m)\right)\eqr q^{rt/2}\,
Q_L(\mu_1,\dots,\mu_m). \label{QL2}
\end{equation}
Here $t$ depends only on the linking matrix $(l_{ij})$ of $L$: $$ t= (r-h) \sum_{1\le
i,j\le m }l_{ij} (g_i|g_j) + 2 \sum_{1\le i,j\le m}l_{ij} (g_i|\mu_j-\rho),$$ with $h$
being the Coxeter number of the Lie algebra $\fg$ (see Table 1).
\end{thm}

For the special cases $\fg=sl_2$ and $\fg =sl_n$, the theorem was proved by Kirby and
Melvin \cite{KM} and  Kohno and Takata \cite{KT1}. In \cite{KM,KT1}, the ``twisting
factor" $q^{rt/2}$ is derived by  direct computations. The proof given in
\cite{integrality} for general simple Lie algebra used a tensor product theorem of
Lusztig.

\begin{cor}\label{symm2}
If $\mu_j-\rho$ is in the root lattice and $\mu_j\in \bCr$, then, $$
Q_L\left(g_1(\mu_1),\dots,g_m(\mu_m)\right)\eqr q^{\frac{r(r-h)}{2}[ \sum l_{ij}
(g_i|g_j)]}\, Q_L(\mu_1,\dots,\mu_m). $$
\end{cor}

\begin{rem} When $(r,d)\not =1$, we can strengthen both symmetry principles using
 the Weyl alcove defined by the {\em long} highest root,  see
\cite{integrality}.
\end{rem}

\subsection{More Integrality}  The result of this subsection is new and will help us to define
quantum invariants of 3-manifold without using the complicated theory of quantum groups
at roots of unity.

Every  finite-dimensional $\UU$-module $V$ decomposes as the direct sum of
$\lambda$-homogeneous components, $\lambda\in X_+$ (recall that we work with general
parameter). Here a $\lambda$-homogeneous component $E_\lambda$ is the maximal submodule
isomorphic to the sum of several copies of $\Lambda_\lambda$. Each $\lambda$-homogeneous
component defines a projection $\pi_\lambda:V\to E_\lambda$ and an inclusion 
$\iota_\lambda : E_\lambda \to V$. 

\begin{figure}[htpb]
$$ \psdraw{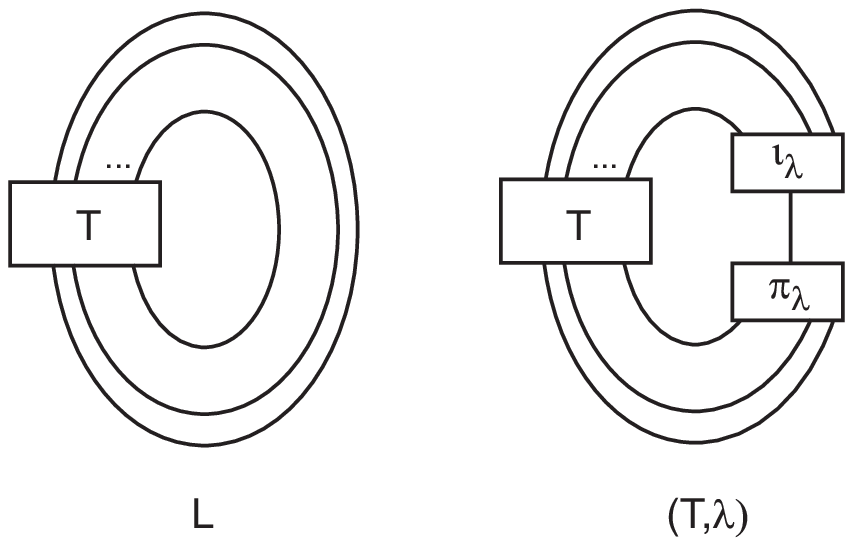}{3 in} $$
\caption{}\label{tangle1}
\end{figure}

Suppose a link $L$ is the closure of a $(n,n)-$tangle $T$, as shown in Figure
\ref{tangle1}a. Let $V_1,\dots,V_n$ are the colors of the strands shown (some of them may
come from the same component). Consider a $\lambda$-homogeneous component $E_\lambda$ of
$V_1\otimes \dots\otimes V_n$. Cut the $n$ strands and insert  2 ``coupons" with 
operators $\pi_\lambda, \iota_\lambda$ in them; we got a ``tangle with coupons"  
$(T,\lambda)$, see Figure \ref{tangle1}b. Ribbon category can be used a define isotopy 
invariant of objects like $(T,\lambda)$, denoted by $J_{(T,\lambda)}$, which is generally 
in the fractional field $\Q(q^{1/2D})$ (see \cite{Turaev}). The following proposition, 
whose proof is based on the result of \cite{integrality} and borrows an idea from 
\cite{MW},  shows that for the special case $(T,\lambda)$, the quantum invariant is a 
Laurent polynomial in $q^{1/2D}$. 

\begin{pro} In the above setting, $J_{(T,\lambda)}$ is in the ring $\Z[q^{\pm 1/2D}]$, and
in that ring, it is
 divisible by the quantum dimension $J_U(\Lambda_\lambda)$.
 \label{s10}
\end{pro}

\begin{pf} According to the general theory, $J_T$ acts as a $\UU$-endomorphism of
$V_1\otimes \dots\otimes V_n$. Hence $J_T$ commutes with $\iota_\lambda\, {}_o \,
\pi_\lambda$. Note that $E_\lambda$ must be  of the form $\Lambda_\lambda\otimes N$, 
where $N$ is a vector space over $\Q(q^{1/2D})$. On $E_\lambda= \Lambda_\lambda\otimes N$ 
the operator  $J_T$ acts as $\id \otimes R$, 
where $R$ is an operator acting on $N$. It follows that

\begin{equation} J_{(L,\lambda)} = \tr(R) \times J_U(\Lambda_\lambda).\label{123}
\end{equation}

Eigenvalues of $R$ are also eigenvalues of $J_T$ which can be represented by a matrix
with entries in $\Z[q^{1/2D}]$ (see \cite{integrality}, the proof used only the theory of 
quantum groups with general parameter). Hence $\tr(R)$ is in the ring ${\cal I}$ of 
algebraic integers over $\Z[q^{1/2D}]$. On the other hand, that the decomposition of 
$V_1\otimes \dots\otimes V_n$ into $\lambda$-homogeneous components can be done over the 
fractional field $\Q(q^{1/2D})$ means  $R$ can be represented by a matrix with entries in 
$\Q(q^{1/2D})$. Thus $\tr(R)\in \Q(q^{1/2D})$. Since ${\cal I}\cap \Q(q^{1/2D})= 
\Z[q^{1/2D}]$, we have that $\tr(R)\in \Z[q^{1/2D}]$. The proposition now follows from 
(\ref{123}). 
\end{pf}

\begin{rem}\label{x32}
We will use the proposition in the following way. First, since $J_{(T,\lambda)}$ is a
Laurent polynomial in $q^{1/2D}$, we can plug  any non-zero value of $q^{1/2D}$ in
$J_{(T,\lambda)}$. Next, for special values of $q^{1/2D}$ annihilating
$J_U(\Lambda_\lambda)$, the value of $J_{(T,\lambda)}$ is 0. Also note that

\begin{equation}
J_L (\mu_1,\dots,\mu_m) = \sum_\lambda J_{(T,\lambda)} \label{sao}
\end{equation}

\end{rem}

\section{Quantum 3-manifold invariants}

\subsection{Introduction}
 Quantum invariants of 3-manifolds can be constructed only when $q$ is a root
of unity of some order $r$. In previously known cases,  $r$  must be divisible by $d$, 
since this will ensure that the so called $S$-matrix is invertible. Our construction of 
3-manifold invariants is slightly in more general situation: we will get invariants of 
3-manifolds even in the case when $r$ is not divisible by $r$. For this reason we will 
give a new (but not quite new) proof of the existence of quantum invariants of 
3-manifolds. We will use only quantum groups with general parameter to define link 
invariants, and  only on the last step we replace $q$ by a root of unity in {\em link 
invariants}. In this paper, a 3-manifold is always closed and oriented.

 For the reader to 
have an idea how the quantum 3-manifold invariant looks like, let us give here the value 
for the Poincare Homology 3-sphere P (obtained by surgery on a left-hand trefoil with 
framing $-1$), with $\fg=sl_2$: 

$$ \tau^{sl_2}_P(q) = \frac{1}{1-q}\sum_{n=0}^\infty q^n (1-q^{n+1})(1-q^{n+2}) \dots
(1-q^{2n+1}).
$$

Here $q$ is  a root of unity, and the sum is easily seen to be finite. Similar, but 
different formula has also been obtained by Lawrence and Zagier, using some calculation 
involving modular forms, see \cite{LZ}. 

When $M$ is the Brieskorn sphere $\Sigma(2,3,7)$ (obtained by surgery on the right hand 
trefoil with framing $-1$): 

$$\tau_M^{sl_2}(q) = \frac{1}{1-q}\sum_{n=0}^\infty q^{-n(n+2)} (1-q^{n+1})(1-q^{n+2})\dots (1-q^{2n+1}).$$
Again $q$ must be a root of unity for the above expression to have meaning.

\subsubsection{Heuristic} The values of $Q_L(\mu_1,\dots,\mu_m)$ are in $\Z[q^{\pm
1/2D}]$. The infinite sum $\displaystyle{\sum_{\mu_j \in X} Q_L(\mu_1,\dots,\mu_m)}$ does
not have any meaning. It is believed  (and there are reasons for this ) that the sum is
invariant under the second Kirby move, and hence almost defines a 3-manifold invariant.
The problem is to regularize the infinite sum $\displaystyle{\sum_{\mu_j \in X}
Q_L(\mu_1,\dots,\mu_m)}$. One solution is based on the fact that at $r$-th roots of
unity, $Q_L(\mu_1,\dots,\mu_m)$ is periodic (the first symmetry principle), so we should 
use the sum with $\mu_j$'s run over the set $\bar C_r$. 

\subsection{Sum over $\bar C_r$}
\subsubsection{General}

Let us fix a positive integer $r \ge d h^\vee$, called the {\em shifted level}, and a
primitive $r$-th root of unity $\xi$. At some stage we also need a primitive $2Dr$-th
root $\zeta$ of $1$ such that $\zeta^{2D}=\xi$. Let

\begin{equation} F^\fg_L(\xi;\zeta) = \sum_{\mu_j\in \bar C_r\cap X} Q_L(\mu_1,\dots,\mu_m)|_{q^{1/2D}=\zeta}.
\notag
\end{equation}
Certainly $\zeta$ determines $\xi$, but we prefer to keep $\xi$ in the notation $
F^\fg_L(\xi;\zeta)$ since in later cases, the whole thing depends only on $\xi$ but not
$\zeta$, and in those cases we will drop $\zeta$ in the notation.

 Recall that $h = (\rho|\alpha_0) +1$ is the Coxeter number. Let $k=r-h$. Then we have
 (recalling the shift by $\rho$ and the fact that $Q_L \eqr 0$ on the boundary of $\bar C_r$)

\begin{equation}
F^\fg_L(\xi;\zeta) = \sum_{\mu_j\in \bar C_k}
Q_L(\Lambda_{\mu_1},\dots,\Lambda_{\mu_m})|_{q^{1/2D}=\zeta}. \label{definitionF2}
\end{equation}

The following half-open parallelepiped $P_r$ is a fundamental domain of the group $rY$:

$$ P_r = \{ x= c_1 \al_1 +\dots + c_\ell \al_\ell \in \hh^*_\R \mid  0\le c_1, \dots,
c_\ell < r\}.$$

Since $\bar C_r$ is a fundamental domain of $W_r= W\ltimes rY$, we have, due to the first
symmetry principle,

\begin{equation} F^\fg_L(\xi;\zeta) = \frac{1}{|W|}\sum_{\mu_j\in \bar P_r\cap X}
Q_L(\mu_1,\dots,\mu_m)|_{q^{1/2D}=\zeta}. \label{definitionF3}
\end{equation}

\subsubsection{The Gauss sum} The following is  a quadratic Gauss sum on the abelian group
$X/rY \equiv P_r\cap X$, with $\zeta$ used to define fractional powers of $\xi$:

$$ \gamma^\fg (\xi,\zeta) := \sum _{\mu \in P_r\cap X} \xi^{\frac{1}{2}(|\mu|^2
-|\rho|^2)}.$$

Criteria for vanishing of a Gauss sum are known, see eg. \cite{Deloup}. Using the
explicit structure of simple Lie algebras and the criterion one can prove the following.

\begin{pro}  The Gauss sum $\gamma^\fg(\xi,\zeta)=0$ if and only if $r$ is odd and $\fg$
is either $C_\ell$ with arbitrary $\ell$,  or $B_\ell$ with even $\ell$. \label{vanish}
\end{pro}

The following lemma uses the well-known trick of completing the square.

\begin{lem}
Suppose $\beta\in X$. Then

$$ \sum_{\mu\in P_r\cap X} \xi^{\frac{1}{2}(|\mu|^2-|\rho|^2)} \xi^{(\beta|\mu)} =
\gamma^\fg(\xi,\zeta) \times \xi^{-\frac{1}{2}|\beta|^2}.$$ \label{square}
\end{lem}

\begin{pf} Completing the square, we see that

$$ \frac{1}{2} (|\mu|^2 -|\rho|^2) + (\beta|\mu) = \frac{}{}(|\mu+\beta|^2 - |\rho|^2)
-\frac{1}{2}|\beta|^2.$$

It remains to notice that everything is invariant under the translation $rY$, and both
$P_r$ and $P_r +\beta$ are fundamental domains of $rY$.
\end{pf}

\subsection{Invariance under the second Kirby move}

\begin{pro}  Suppose that the order $r$ of $\xi$ is greater than or equal to $d h^\vee$. Then
 $F^\fg_L(\xi,\zeta)$ does not depend on the orientation of $L$ and is invariant under the
 second Kirby move.
 \label{inva}
\end{pro}

\begin{pf}
Using linearity we extend the invariant $J_L$ to the case when the colors are elements of
the $\Z[q^{\pm 1/2D}]$-module freely generated by $\Lambda_\lambda,\lambda\in X_+$. Then

$$ F^\fg_L(\xi,\zeta) = J_L(\omega,\dots,\omega)|_{q^{1/2D}= \zeta},$$ where

$$ \omega = \sum_{\mu \in  int(C_r) \cap X} J_U(\mu)\, \Lambda_{\mu-\rho} =
 \sum_{\mu \in  \bar C_k \cap X} J_U(\Lambda_\mu)\, \Lambda_{\mu} .$$

The independence of orientation  is simple: If we reverse the orientation of one
component, and at the same time change the color from $V$ to the dual $V^*$, then the
quantum link invariant remains the same. It is known that the alcove $\bar C_k$, (here
$k=r-h$), is invariant under taking dual, i.e. the dual of $\Lambda_\mu, \mu \in \bar
C_k$ is another $\Lambda_{\mu^*}$, with $\mu^*$ again in $\bar C_k$. Moreover,
$J_U(\Lambda_\mu)= J_U(\Lambda_\mu^*)$. Hence $\omega$ is invariant under $\mu\to \mu^*$,
and $J_L(\omega,\dots,\omega)$ is unchanged if we reverse the orientation of one
component.

Let us consider the 2-nd Kirby move $L \to L'$, as described in Figure \ref{Kirby2}, with
blackboard framing.  In both $L,L'$ let $K$ be the singled out unknot component with 
framing 1. 
\begin{figure}[htpb]
$$ \psdraw{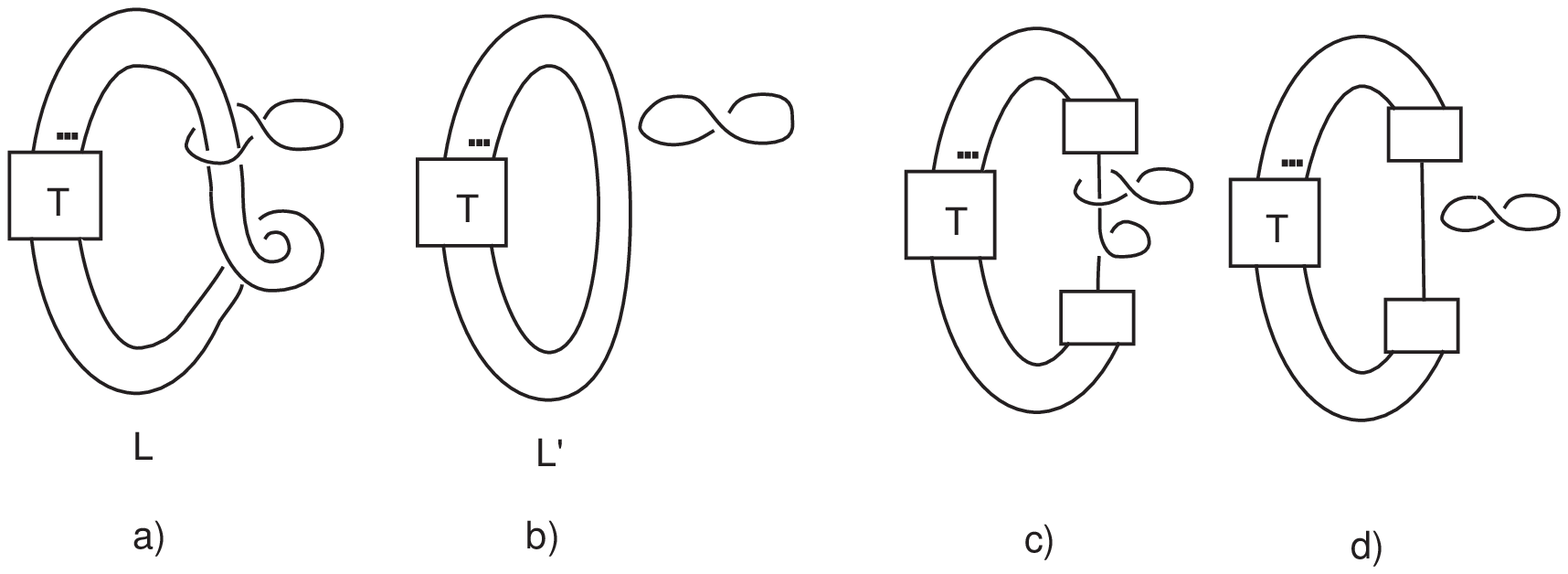}{4.6in} $$
\caption{}\label{Kirby2}
\end{figure}

Then we have to show that

$$ J_L(\omega,\dots,\omega) \eqr J_{L'}(\omega,\dots,\omega).$$

It is enough to show that for every $\mu_1\dots,\mu_m\in X_+$,

$$J_L(\mu_1,\dots,\mu_m, \omega) \eqr J_{L'}(\mu_1,\dots,\mu_m,\omega).$$
 Here we suppose $L$ and $L'$ have $m+1$ components with $K$ being the $m+1$-st.

Suppose the colors of the $n$ strands coming out from the box $T$ are $V_1, \dots,V_n$.
(Each $V_i$ is one of $\Lambda_{\mu_j-\rho}$ or their duals.) The module
$V_1\otimes\dots\otimes  V_n$ is completely reducible over $\Q(q^{1/2D})$, so we
decompose it into homogeneous components. Using (\ref{sao}) to decompose $J_L$ and 
$J_{L'}$ into sums of quantum invariants of ``tangles with coupons", see Figure 
\ref{Kirby2}c,d. In each tangle with coupons there is only one strand, with color a 
homogeneous component, piercing through $K$. Now put $q^{1/2D}=\zeta$ (see Remark 
\ref{x32}). We see that it's remain to prove the following lemma, which is essentially 
the statement of the proposition for the case when $n=1$. 

\begin{lem}
Suppose $J_U(\lambda)|_{q^{1/2D} =\zeta} \neq 0$. Then

$$ J_Z(\lambda,\omega)|_{q^{1/2D} =\zeta}  =  J_{Z'} (\lambda,\omega)|_{q^{1/2D}
=\zeta},$$ where $Z,Z'$ are the $(1,1)$-tangles in Figure \ref{tangle2}.
\end{lem}

\begin{figure}[htpb]
$$ \psdraw{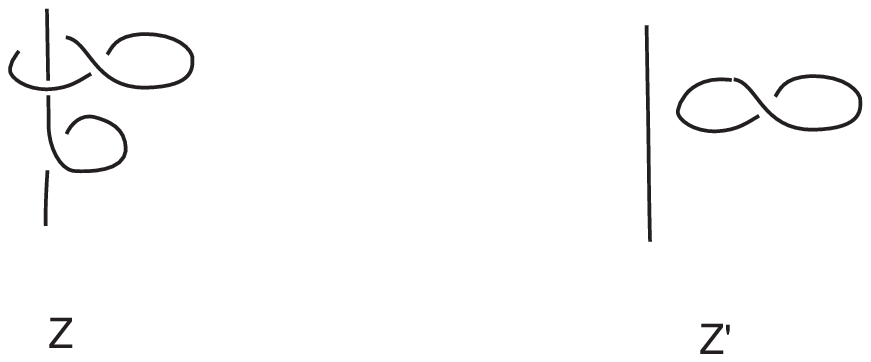}{2.6in} $$
\caption{}\label{tangle2}
\end{figure}

Proof.  Note that both sides of $J_Z(\lambda,\omega)$ and $J_{Z'}(\lambda,\omega)$ are
scalar operator acting on $\Lambda_{\lambda-\rho}$. Closing $Z$, we get the Hopf link
$H_+$ with framing 1 on both components. Similarly, closing $Z'$ we get the trivial link
$U_2$ with framing 1 on the second component. Since $J_U(\lambda)|_{q^{1/2D}=\zeta} \neq 
0$, the identity to prove is equivalent to 

\begin{equation} J_{H_+} (\lambda,\omega) \eqr J_{U_2}(\lambda,\omega).
\label{sao1}
\end{equation}

Let us first calculate

$$ u = \sum_{\mu \in P_r \cap X}  q^{\frac{|\mu|^2 -|\rho|^2}{2}}J_U (\mu)
\,J_H(\lambda,\mu)|_{q^{1/2D} =\zeta }.$$

We have, with $q^{1/2D} = \zeta$,

\begin{align} u
&=
\frac{1}{\psi^2} \sum_{\mu\in P_r \cap X} \sum_{w,w'\in W}\sn(w w')\,
q^{\frac{|\mu|^2-|\rho|^2}{2}} q^{(\mu|w(\lambda)+ w'(\rho))} \qquad \text{by formula
(\ref{Hopflink})}\notag \\
&=  \frac{1}{\psi^2}
\gamma^\fg(\xi,\zeta) \sum_{w,w'\in W } \sn (w w') \,q^{-\frac{1}{2}(|\lambda|^2 +
|\rho|^2 + 2 (\lambda|w^{-1} w'(\rho) )} \qquad \text{by Lemma \ref{square}}\notag\\
&= \frac{1}{\psi^2} \gamma ^\fg (\xi,\zeta)\, q^{-\frac{
|\rho|^2 + |\lambda|^2}{2}} |W| \sum_{w\in W} \sn (w) q^{(-\lambda|w(\rho))}\notag\\
&=   \frac{1}{\psi}   q^{-\frac{
|\rho|^2 + |\lambda|^2}{2}} |W|\, \gamma ^\fg (\xi,\zeta)\, J_U(-\lambda) \qquad \text{by
(\ref{unknot2})} \notag \end{align}

Thus

\begin{equation}
\sum_{\mu \in P_r \cap X}  q^{\frac{|\mu|^2 -|\rho|^2}{2}}J_U (\mu) \,J_H(\lambda,\mu) 
\eqr \frac{(-1)^s |W|}{\psi}   q^{-\frac{ |\rho|^2 + |\lambda|^2}{2}} |W| \gamma ^\fg 
(\xi,\zeta)\, J_U(\lambda).
\label{Smatrix}
\end{equation}

 Recall that for $Q_L$,  increasing  by 1 the framing of a component
colored by $\Lambda_{\lambda -\rho}$ results in a factor $q^{(|\mu|^2 -|\rho|^2)/2}$. The 
left hand side of (\ref{sao1}) is 

$$ LHS = \frac{1}{|W|} q^{\frac{|\mu|^2-|\rho|^2}{2}} \sum_{\mu \in P_r \cap X}  q^{\frac{|\mu|^2 -|\rho|^2}{2}}J_U (\mu)
\,J_H(\lambda,\mu)\mid_{q^{1/2D} =\zeta }.$$

The right hand side is

\begin{align} RHS &= J_U(\lambda) \frac{1}{|W|} \sum_{\mu \in P_r \cap X} q^{\frac{|\mu|^2-|\rho|^2}{2}}
J_U(\mu)^2\mid_{q^{1/2D}=\zeta} \notag \\
&=  J_U(\lambda) \frac{1}{|W|} \sum_{\mu \in P_r \cap X} q^{\frac{|\mu|^2-|\rho|^2}{2}} J_H(\mu,\rho)
\, J_U(\mu)\mid_{q^{1/2D}=\zeta}.\notag
\end{align}

Hence it follows from (\ref{Smatrix}) that $LHS =RHS$.  \end{pf}

Let  us record here the formula for  $F_{U_+}(\xi;\zeta)$, where $U_+$ is the unknot with
framing 1.

$$ F^\fg_{U_+}(\xi;\zeta) = \frac{1}{|W|} \sum_{\mu\in P_r \cap X}
 q^{\frac{|\mu|^2-|\rho|^2}{2}} J_H(\rho,\mu) J_U(\mu)\mid_{q^{1/2D}=\zeta}
,$$ and hence (\ref{Smatrix}) gives

\begin{equation}
F^\fg_{U_+}(\xi;\zeta) =  \frac{\gamma^\fg(\xi;\zeta)}{\prod_{\al
 >0}(1-q^{{(\al|\rho)}})}.
 \label{FU}
 \end{equation}

\begin{rem}  In the proof we used the first symmetry principle, whose proof required the 
theory of quantum groups at roots of unity. However, if we defined $F_L$ using the sum 
over $P_r \cap X$ at  the beginning, then we would not have to use the first symmetry 
principle. 
\end{rem} 

\subsection{Quantum Invariants}
\subsubsection{Definition}
 Suppose $U_{\pm}$ are the unknot with framing $\pm1$. Note that
$F^\fg_{U_{\pm}}(\xi;\zeta)$ are complex conjugate to each other. If
$F^\fg_{U_{\pm}}(\xi;\zeta)\neq 0$, then one can define invariant of the 3-manifold $M$
obtained by surgery along $L$ by the formula:

$$\tau^\fg_M(\xi;\zeta) := \frac{F_L(\xi;\zeta)}{F_{U_{+}}(\xi;\zeta)^{\sigma_+}\,
F_{U_-}(\xi;\zeta)^{\sigma_-}}.$$  Here $\sigma_+,\sigma_-$ are the number of positive
and negative eigenvalues of the linking matrix of $L$. If $F^\fg_{U_{\pm}}(\xi;\zeta)=0$,
then let $\tau_M^\fg(\xi;\zeta)=0$ for every 3-manifold $M$.

Here are the cases when $F^\fg_{U_{\pm}}(\xi;\zeta) =0$.

\begin{pro}  Suppose the order $r$ of the root $\xi$ satisfies $r \ge d h^\vee$. Then
$F^\fg_{U_{\pm}}(\xi;\zeta)=0$ if and only $r$ is odd and $\fg$ is either $B_\ell$ with
even $\ell$ or $C_\ell$ with arbitrary $\ell$.  In particular, if $r$ is divisible by
$d$, then $F^\fg_{U_{\pm}}(\xi;\zeta)$ is not equal to 0.
\label{not0}
\end{pro}

\begin{pf}  The proposition
follows from formula (\ref{FU}) and Proposition \ref{vanish}. \end{pf}

\begin{rem} Only in the two cases listed in the proposition are the invariants trivial. But $F_{U_\pm}\neq 0$
does mean that the so-called $S$-matrix is invertible.
\end{rem}

\subsubsection{Comparison with known cases} In the literature, only the case $r$ divisible
by $d$ was considered. In that case $r = d r'$, and the number $r'-h^\vee$ is called the
level of the theory (see \cite{Kirillov}). Also in this case one can construct a modular 
category, and a  topological quantum field theory. 

Here we  consider both cases when $r$ is or is not divisible by $d$. In the latter case,
the level should be $r-h$.

 In the book \cite{Turaev}  modular category, and hence
quantum invariants, was constructed for  simple Lie algebras with $d=1$. Later work of
\cite{AP} established the existence of modular category for every simple Lie algebra, at
shifted level $r$ divisible by $d$, see a rigorous proof in \cite{Kirillov}. We will 
explain here why the invariant of \cite{Kirillov} is coincident with ours, when $r$ is 
divisible by $d$. 

If $d=1$, then the set of modules $\Lambda_\mu$, with $\mu+\rho\in C_r$  forms a modular 
category (see \cite{AP,Kirillov}), hence the 3-manifold invariant derived from the 
modular category is exactly our  $\tau^\fg_M(\xi;\zeta)$.  

Suppose $d>1$, and $r$ is divisible by $d$. In this case the above set of modules does 
not form a modular category. There is a smaller simplex $C'_r \subset C_r$ with the 
corresponding affine Weyl group $W'_r$ such that $C_r$ consists of several copies of 
$C'_r$ under the action of $W'_r$; and the modular category  consists of $\Lambda_\mu$ 
with $\mu+\rho \in C'_r$. The corresponding 3-manifold invariant is thus obtained by 
taking the sum over the smaller simplex $C'_r$. The first symmetry principle is valid if 
$C_r,W_r$ are replaced with $C'_r, W'_r$ (for details see \cite{integrality}). 
 Due to this symmetry, the sum of
$Q_L$ over the bigger simplex $C_r$ is simply a constant  times the sum over $C'_r$. This
is the reason why we can use $C_r$ to define the same 3-manifold invariant. This smaller 
simplex $C'_r$ is constructed using the {\em long} highest root.

\section{Quantum invariant of the projective group}

\subsection{Preliminaries}
There is a simply-connected complex Lie group $\cal G$ corresponding to $\fg$. The
invariant $\tau^{\fg}_M$ is associated  to $\cal G$. Let $G$ be the center group of $\cal
G$. It is known that $G$ is isomorphic to $X/Y$, and $|G|=\det(a_{ij})$. For every
subgroup $G'\subset G$, there corresponds a Lie group ${\cal G}/G'$, and there is a
quantum invariant associated with this quotient group. We will describe here a method to
construct them, focusing on the extreme case when $G'=G$.  We will see that there are
many shifted levels $r$ for which the invariant $\tau^\fg_M$ is trivial, but at the same 
time the invariant of the projective group, denoted by $\tau^{P\fg}_M$, is non-trivial, 
and even defined by a modular category. We will see that if $r$ and $\det(a_{ij})$ are 
coprime, i.e. $(r,\det(a_{ij}))=1$, then the invariant associated to the projective group 
is not trivial. 

\subsubsection{The lattice $\rho +Y$}
\begin{lem} For every positive integer $r$, the lattice $\rho+Y$ is invariant under
the action of $W_r$. \label{lattice}
\end{lem}
\begin{pf} Recall that $W_r = W\ltimes rY$. The fact that $\rho+Y$ is invariant under
$rY$ is obvious. That $\rho+ Y$ is invariant under the action of $W$ follows from the
fact that $w(\rho)-\rho$ belongs to $Y$. (Actually, $w(\mu)-\mu \in Y$ for every $\mu\in
X$.)
\end{pf}

\subsubsection{Sums over the root lattice}

Let $k=r-h$, and $\xi$ is a root of unity of order $r \ge d h^\vee \ge h$. Let

$$ F^{P\fg}_L(\xi) = \sum_{\mu_j\in (\bar C_k\cap Y)}
Q_L(\Lambda_{\mu_1},\dots,\Lambda_{\mu_m})|_{q=\xi}.$$

The definition is the same as in (\ref{definitionF2}), except that we sum over $\mu_j$'s
which are in the {\em root lattice}. Note that there is no need to fix a $2D$-th root of 
$\xi$, since by Corollary \ref{integ}, there is no fractional power of $q$. 

Recalling the shift by $\rho$, we have

$$ F^{P\fg}_L(\xi) = \sum_{\mu_j\in \bar C_r\cap
(\rho+Y)} Q_L({\mu_1},\dots,{\mu_m})|_{q=\xi}.$$

Lemma \ref{lattice} and the first symmetry principle show that

 $$ F^{P\fg}_L(\xi) =  \frac{1}{|W|}\sum_{\mu_j\in \rho +(P_r\cap Y)}
Q_L({\mu_1},\dots,{\mu_m})|_{q=\xi}.$$

\subsubsection{Gauss sum}
We will encounter a Gauss sum on the group $Y/rY$.  From now on  let  
$${\sum}_r$$ stands
for $\displaystyle{\sum_{\mu\in \rho +( P_r \cap Y)}}$. Put

$$ \gamma^{P\fg}(\xi) := {\sum}_r\,
\xi^{\frac{|\mu|^2-|\rho|^2}{2}}.$$

Then, the same proof of  Lemma \ref{square} gives us:

\begin{lem}  Suppose $\beta \in Y$. Then

$$   {\sum}_r\, \xi^{\frac{|\mu|^2-|\rho|^2}{2}} \xi^{(\mu|\beta)} =
\gamma^{P\fg}(\xi) \times \xi^{-\frac{|\beta|^2}{2}}.$$
\label{square2}
\end{lem}

\subsection{Definition of invariants associated to the projective group}
 Recall that $H$ is the Hopf link. Let $S_{\lambda,\mu} = J_H(\lambda,\mu)|_{q=\xi}$.
Let the matrix $S$ have entries $S_{\lambda,\mu}$ with $\lambda,\mu\in 
\text{Interior}(C_r) \cap (\rho+Y)$. 

\begin{thm}\label{212x}

a) Suppose the order $r$ of $\xi$ is greater than or equal to $d h^\vee$. Then
$F^{P\fg}_L(\xi)$ is invariant under the 2-nd Kirby move and does not depend on the
orientation of $L$.

b) If $r$ is coprime with $d \det (a_{ij})$, then the matrix $S$ is invertible.

c)  If $r$ is coprime with $\det (a_{ij})$, then $F_{U_\pm}^{P\fg}(\xi) \neq 0$.

\end{thm}

\begin{pf} Notice that if $\lambda,\mu\in Y$, then in the decompostion of
$\Lambda_\lambda\otimes \Lambda_\mu$ into irreducible modules one encounters only
$\Lambda_\nu$ with $\nu\in Y$. This is a well-known fact: The irreducible modules of the
group ${\cal G}/G$ have highest weights in $X_+\cap Y$, and finite-dimensional ${\cal
G}/G$-modules are completely reducible.

Using this fact one can repeat the proof of Proposition \ref{inva} to get a proof of part
a).

b) We will show that $S \bar S$ is a non-zero constant times the identity matrix. Here
$\bar S$ is the complex conjugate. We know that $\psi \neq 0$ when $q=\xi$, since $r \ge
dh^\vee$ (see \ref{trivial}).  Using (\ref{Hopflink}) and $\sum_{\mu\in \bar C_r\cap 
(\rho +Y) }=\frac{1}{|W|}{\sum}_r$, we have 

\begin{align} {|W|}\psi ^2 (S \bar S)_{\lambda,\nu}&=
 {\sum}_r\,\sum_{w,w' \in W} \sn(w w')
 \xi^{(\mu|w(\lambda)-w'(\nu))}\notag \\
 & = \sum_{w,w' \in W} \sn(w w') \left[ {\sum}_r\,
 \xi^{(\mu|w(\lambda)-w'(\nu))}\right].\notag \end{align}

 Let $Y^*$ be the lattice dual  to $Y$, over $\Z$,
 with respect to the scalar product. If  $w(\lambda)-w'(\nu) \not \in rY^*$,
 then there is a fundamental root $\al_i$ such that $(\al_i|w(\lambda)-w'(\nu)) \not \in r\Z$. It follows
 that the sum in the square bracket is 0, since $\xi^n + \xi^{2n} + \dots
 +\xi^{(r-1)n}=0$ if $n$ is not divisible by $r$, and ${\sum}_r$ is the sum over a fundamental domain of $rY$.

We will find out when $w(\lambda)-w'(\nu) \in rY^*$. Note first that $w(\lambda)
 -w'(\nu) \in Y$. We'll find the intersection $r Y^*\cap Y$.

The lattice $Y^*$ is spanned by $\lambda_1/d_1,\dots,\lambda_\ell/d_\ell$, where
$\lambda_1,\dots,\lambda_\ell$ are the fundamental weights. Thus the order of $Y^*/X$ is
$d_1 d_2 \dots d_\ell$, a factor of $d^\ell$. The order of $X/Y$ is $\det (a_{ij})$. Thus
the group $rY^*/rY \equiv Y^*/Y$ has order a factor of $d^\ell \times \det (a_{ij})$.

The group $Y/rY$ has order $r^\ell$. By
 assumption, the orders of two groups  $r Y^*/rY$ and $Y/rY$ are co-prime. Their intersection must be
 trivial.  Hence $w(\lambda)-w'(\nu)$, belonging to both $r Y^*$ and $Y$, must belong
 to $rY$.  But this means $\lambda$ and $\nu$ are in the same $W_r$-orbit. This could
 happen for $\lambda,\nu \in Int(C_r)$ if and only if $w=w'$ and $\lambda =\nu$.

 When $w=w'$ and $\lambda=\nu$, the sum in the square bracket is $r^\ell$.
Thus $(S \bar S)$ is a non-zero constant times the identity.

c) One can  prove c) directly using the
 criterion of vanishing of Gauss sums. Or one can use the following arguments. If, in
 addition, $r$ is coprime with $d$, then by b), the $S$ matrix is non-degenerate. In this
 case it is known that $F_{U_\pm}^{P\fg}(\xi) \neq 0$ (see \cite{Turaev}). Suppose now
 $(r,d)\neq 1$. Then $r$ is divisible by $d$.
Formula (\ref{sao5}) below shows that $F_{U_+}^{P\fg}(\xi)$ is a
 factor of $F^\fg_{U_+}(\xi,\zeta)$, which is not 0 by Proposition
 \ref{not0} (for some $2D$-th root $\zeta$ of $\xi$).  Hence  $F_{U_\pm}^{P\fg}(\xi) \neq 0$.
\end{pf}

If $F^{P\fg}_{U_\pm}(\xi) =0$, we define $\tau_M^{P\fg}(\xi)=0$, otherwise, let

\begin{equation}
\tau^{P\fg}_M(\xi) := \frac{F^{P\fg}_L(\xi)}{(F^{P\fg}_{U_{+}}(\xi))^{\sigma_+}
(F^{P\fg}_{U_-}(\xi))^{\sigma_-}}, \label{definition3}\end{equation}
 where $M$ is
obtained from $S^3$ by surgery along the framed link $L$.

\begin{rem} 

a)  Theorem \ref{212x}, part (b) shows that when $r$ is co-prime with $d \det(a_{ij})$, 
the set of all modules $\Lambda_\mu$, with $\mu\in \bar C_k\cap Y$ (note the root lattice 
$Y$ here), generates a modular category. Here one has  to use the reduced quotient 
structure as in \cite{AP,Kirillov}. At the same time, if $Y$ is replaced by $X$, then the 
resulting category, usually considered by algebraists (says, in earlier papers of H. 
Andersen) might not be a modular category. The reason  is the $S$-matrix might not be 
invertible. There are values of $r$ when the $S$-matrix  is invertible for the $Y$ case, 
but  not for the $X$ case. 

b) Whenever $F_{U_\pm}^{P\fg}(\xi) \neq 0$, one has non-trivial invariants. In addition 
to the cases described in the theorem, there are other cases when $F_{U_\pm}^{P\fg}(\xi) 
\neq 0$. 
 For example, using the 
criterion for the vanishing of Gauss sum, one can also prove that whenever $r$ is odd 
(for all $\fg$), $F_{U_\pm}^{P\fg}(\xi) \neq 0$. On the other hand, there are cases when 
$F_{U_\pm}^{P\fg}(\xi) = 0$: Examples include the case $\fg =sl_2$, $r$ is divisible by 
4.

b) The invariant $\tau_M^{P\fg}$ is the invariant associated with the projective group,
since the root lattice is spanned by highest weights of  finite-dimensional irreducible
modules of ${\cal G}/G$.  If $G'$ is a subgroup of $G$, then one can construct invariant
associated with ${\cal G}/G'$ by using the lattice $Y'$ generated by the set of all
highest weights of ${\cal G}/G'$. The construction is similar.

\end{rem}

\subsection{Invariant associated to a finite abelian group with a bilinear form} On the
group $G=X/Y$ there is defined the symmetric bilinear form $(\cdot|\cdot)$ with values in
$\frac{1}{D}\Z/\Z \subset \Q/\Z$. For any such group there is a way to define invariants
of 3-manifolds which carry only the information about the homology groups and the linking
form on the torsion of the first homology group, see \cite{MOO,Deloup,Turaev}. We will
present here the theory in the form most convenient for us.

Again $\zeta$ is root of unity of order $2Dr$, and $\xi= \zeta^{2D}$. Define

$$ F^G_L(\xi;\zeta) := \sum_{g_i,g_j\in G} \xi^{r(r-h) \times \frac{1}{2}\sum l_{ij}
(g_i|g_j)},$$ where $l_{ij}$ is the linking matrix of $L$. Here we use $\zeta$ to define
fractional powers of $\xi$.

Then $$F^G_{U_\pm}(\xi;\zeta) = \sum_{g\in G} \xi ^{ r(r-h) (g|g)/2}$$ is a Gauss sum. If
$F^G_{U_\pm}(\xi;\zeta) =0$, we define $\tau^G_M(\xi;\zeta)=0$, otherwise we define, for

$$ \tau^G_M(\xi;\zeta) := \frac{F^G_L(\xi;\zeta) }  { (F^G_{U_+})^{\sigma_+}\,
(F^G_{U_-})^{\sigma_-}  },$$

for $M$ obtained by surgery on a framed link $L$. It is a 3-manifold invariants. In 
general, $\zeta^{r(r-h)}$ is a root of unity of order $2D$. If $\zeta^{r(r-h)}= \exp(2\pi 
i/2D)$, then our invariant is coincident with those in \cite{Deloup}. 

\subsection{Splitting}

\begin{lem} Suppose $(r,\det (a_{ij}))=1$. Then

a) $G$ acts freely on the set $\bar C_r \cap X$.

b) In each $G$-orbit of $\bar C_r \cap X$ there is exactly one element in $\rho + Y$.

\end{lem}

\begin{pf}

a) Note that $\bar C_r \cap X$ is a finite set. Suppose  $g(\mu)=\mu$ for some $\mu \in
\bar C_r \cap X$, we will show that $g$ is the identity of $G$. There is a lift $\tilde
g\in X$ of $G$ such that

$$ r\tilde g + \mu = \mu \pmod{W_r},$$
which, due to $W_r= W\ltimes rY$, means there is $w\in W$ such that

$$ r\tilde g + \mu \in  w(\mu) + rY.$$

Since $w(\mu)-\mu \in Y$, it follows that $r \tilde g\in Y$, or $r g=0$ in $G$. Because
$(r,|G|)=1$, this implies $g=0$ in $G$.

b) Using $X/Y = rX/Y$ (since $(r,|G|)=1$), we have

$$ X =  Y+ rX,$$

Hence $(\rho+ Y) + rX = X$.  This shows that in each $G$-orbit there is at least one
element in $\rho +Y$. The proof of part a) shows that each $G$-orbit contains at most one
element in $\rho+Y$. \end{pf}

Suppose $(r, \det (a_{ij}))=1$. By the above lemma and the second symmetry principle (see
Corollary \ref{symm2}), one has

\begin{equation} F^\fg _L(\xi;\zeta) = F_L^G(\xi;\zeta) F_L^{P\fg}(\xi).
\label{sao5}
\end{equation}

Hence we have the following splitting theorem

\begin{thm} Suppose $(r, \det (a_{ij}))=1$ and $\zeta$ is a $2Dr$-th root of unity, $\xi=
\zeta^{2D}$. Then $$ \tau_M^\fg (\xi;\zeta) = \tau^{P\fg}_M(\xi) \,
\tau^G_M(\xi;\zeta).$$
\label{split}
\end{thm}

\begin{rem}

a) The invariant  $\tau^G_M(\xi;\zeta)$ carries only the information about the first
homology group and the linking form on its torsion; it is a weak invariant, and sometimes
it is equal to 0, in which case $\tau_M^\fg (\xi;\zeta)=0$. Hence $\tau^{P\fg}_M(\xi)$ is
in general a finer invariant. For example, if $\fg = B_\ell$, and $r$ is odd, then
$\tau_M^\fg (\xi;\zeta)=0$, but $\tau^{P\fg}_M(\xi)$ is in general not 0.

b) When $r$ is not coprime with $\det (a_{ij})$, there are cases when both $\tau^\fg$ and
$\tau^{P\fg}$ are non-trivial, but there is no simple relation between the two
invariants. Examples of such case are: $\fg = sl_n$ and $(r,n) \neq 1$, $\fg =D_\ell$ and 
$r$ even, and $\fg =C_\ell$  with $\ell$ odd and $r$ even. 

c)  The splitting of Theorem \ref{split} fits very well with the Gussarov-Habiro theory
of finite type 3-manifold invariants: In that theory one has first to partition the set
of 3-manifolds into subset of ones with the same homology and linking form, then defines
finite type invariants in each subset using a suitable filtration. The invariant $\tau^G$
corresponds to homology and the linking form, and $\tau^{P\fg}$ can be expanded into
power series, at least for rational homology 3-spheres (see below), that gives rise to
finite type invariants.

 d)  The projective quantum invariants were defined and the splitting theorem was proved
in Kirby-Melvin \cite{KM} for $\fg =sl_2$ and Kohno-Takata \cite{KT1} for $\fg=sl_n$. For
the case when $r$ is divisible by $d$,  a similar splitting has also been obtained by
Sawin \cite{Sawin}, but his proof does not go through for all simple Lie algebras, he has
to exclude a half of $D$ series. For $d=2$, Sawin's result and Theorem \ref{split}
address {\em different} cases, and hence they complement each other.
\end{rem}

\section{Integrality}
\begin{thm} Suppose that $r\ge d h^\vee$ is a prime and   not a factor of
$|W|\det (a_{ij})$, and $\xi$ a primitive $r$-th root of unity. Then $\tau_M^{P\fg}(\xi)$
is in $\Z[\xi]=\Z[\exp(2\pi i/r)]$.
\label{integ3}
\end{thm}

 The theorem was proved in the $\fg =sl_2$ case by H. Murakami \cite{Murakami}
(see also \cite{MR}) and $\fg =sl_n$ by Takata-Yokota \cite{TY} and Masbaum-Wenzl
\cite{MW}. It is conjectured that even when $r$ is a not prime, one also has
$\tau_M^{P\fg}(\xi)\in \Z[\xi]$. The remaining part of this section is devoted to a proof
of this theorem.

\subsection{General facts}
For $a, b \in \Z[\xi]$, we write $a\sim b$ if there is a unit $u$ in $\Z[\xi]$ such that
$a = u b$. Suppose $r$ is an odd prime. It is known that    $(\xi-1)$ is prime in
$\Z[\xi]$, and $r \sim (\xi-1)^{r-1}$. It follows that $(r-1)!$ is coprime with
$(\xi-1)$. If $(n,r)=1$ then  $(\xi^n -1) \sim (\xi-1)$.

Formula (\ref{unknot1}) shows that for every $\lambda \in Y$,
  $$J_U(\Lambda_\lambda)|_{q= \xi}
\sim 1.$$

To prove the theorem, we have to show that the numerator of the right hand side  of
(\ref{definition3}) is divisible by the denominator.  First we will show that the
denominator is just a power of $(\xi-1)$, then we show that the numerator can be
decomposed as a sum of simple terms, each divisible by that power of $(\xi-1)$.

\subsection{Gauss Sum again} Suppose $b$ is an integer. Recall that ${\sum}_r$ stands for
$\sum_{\mu \in \rho+(P_r \cap Y)}$. Let

$$ \gamma_{b}^{P\fg}(\xi) = {\sum}_r \, \xi^{b \frac{|\mu|^2 - |\rho^2|}{2}}.$$

Note that for $\mu\in \rho+ Y$, $|\mu|^2-|\rho|^2$ is always an even number. Hence
$\gamma_{b}^{P\fg}(\xi) \in \Z[\xi]$.

\begin{lem} \label{Gammab}
Suppose $r$ is an odd prime not a factor of $ d \det (a_{ij})$. Then
$\gamma_{b}^{P\fg}(\xi)$ is divisible by $(\xi-1)^{ \frac{r-1}{2}\ell}$. Moreover, if $b$
is not divisible by $r$, then $$ \gamma_{b}^{P\fg}(\xi) \sim (\xi-1)^{
\frac{r-1}{2}\ell}.$$
\end{lem}

\begin{pf}
If $b$ is divisible by $r$, then $\gamma_{b}^{P\fg}(\xi) = r^\ell \sim
(\xi-1)^{\ell(r-1)}$, and we are done.

Suppose $b$ is not divisible by $r$. Then $\xi^b$ is  a root of 1 of order $r$. Hence
there is a Galois automorphism $\sigma$ of the field $\Z(\xi)$ over $\Q$ such that
$\sigma(\xi^b) = \exp(2\pi i /r)$.  Since $\sigma(\xi-1) \sim \xi-1$, it's enough to
prove the lemma in the case $\xi^b = \exp(2\pi i/r)$. In this case

$$ \gamma_{b}^{P\fg}(\xi) = \sum_{\mu \in P_r\cap Y} \exp[\frac{\pi i}{r}
(|\mu+\rho|^2-|\rho|^2)].$$

Since $r$ is odd, and $(\mu|\rho)\in \Z$, one has

$$ |\mu+\rho|^2 -|\rho|^2 \equiv |\mu + (r+1) \rho|^2 - (r+1)^2|\rho|^2 \pmod {2r}.$$

It follows that

$$ \gamma_{b}^{P\fg}(\xi) = \exp[\frac{-\pi i}{r} (r+1)^2|\rho|^2]\sum_{\mu \in P_r\cap Y} \exp[\frac{\pi i}{r} ( |\mu + (r+1)
\rho|^2   )].$$

Notice that $(r+1)\rho\in Y$ since $r+1$ is even, and use the $rY$-invariance, we have

$$ \gamma_{b}^{P\fg}(\xi) = \exp[-\frac{\pi i}{r} (r+1)^2 |\rho|^2 ] \times \sum_{\mu \in
P_r\cap Y} \exp \frac{\pi i |\mu|^2}{r}.$$

The first factor is a unit in $\Z[\xi]$. If $P$ is the matrix $(d_i a_{ij})$  (so that
$(\al_i|\al_j) = P_{ij}$), then the second factor is

$$  \sum_{\mu \in P_r\cap Y} \exp \frac{\pi i |\mu|^2}{r} = \sum _{\vec k \in
(\Z/r\Z)^\ell} \exp[\frac{\pi i}{r} \vec k^t P \vec k].$$

It is known that this Gauss sum is $\sim (\xi-1)^{ \frac{r-1}{2}\ell}$. (This  fact can
be proved by diagonalizing the matrix $P$ and use the value of the 1-variable Gauss sum.
The matrix $P$ is non-degenerate over $\Z/r\Z$ since $\det P$ and $r$ are coprime.)
\end{pf}

\begin{lem}  a) Suppose $b$ is an integer coprime with $r$. Let $b^*$ be an integer such
that  $b b^* \equiv 1\pmod r$. Then

\begin{equation}
 {\sum}_r \xi^{b \frac{|\mu|^2 - |\rho^2|}{2}} \, \xi^{(\mu|\beta)}= \xi^{-b^*
\beta^2/2}\, \gamma_{b}^{P\fg}(\xi), \label{quadratic}
\end{equation}

b)  Suppose $r$ is an odd prime, then the left hand side of (\ref{quadratic}) is
divisible by $(\xi-1)^{ \frac{r-1}{2}\ell}$.
\end{lem}

\begin{pf} a)  The proof is similar to that of Lemma \ref{square}, using the trick of completing
the square.

b) If $b$ is not divisible by $r$, then the statement follows from part a) and Lemma
\ref{Gammab}. Suppose $b$ is divisible by $r$. Then the LHS is either 0 or $r^\ell$,
which is $\sim (\xi-1)^{\ell(r-1)}$.
\end{pf}

\subsection{Unknots and simple lens spaces}

 Let $U_b$ be the unknot with framing $b$. We will first find the prime factors of
 $F^{P\fg}_{U_b}$.
\begin{pro}  a) Suppose $r\ge d h^\vee$ and is coprime with $b$. Let $b^*$
be an integer such that $b b^* \equiv 1 \pmod r$.
 Then

\begin{equation}
 F^{P\fg}_{U_b}(\xi) =  \frac{ \xi^{(1-b^*) |\rho|^2}\, \gamma^{P\fg}_{b}(\xi)\,
J_U(b^*\rho)}{\prod_{\al>0} (1 -\xi^{(\al|\rho)})} .
\label{x2}
\end{equation}

b)   If, in addition, $r$ is an odd prime, then $F^{P\fg}_{U_b}(\xi) \sim (\xi-1)^{(r
\ell -dim \fg)/2}$.

\label{ss23}
\end{pro}

\begin{pf}
a) The proof is similar to that of (\ref{FU}): with $q=\xi$ in $\psi$,

\begin{align}
 F^{P\fg}_{U_b}(\xi) &= \frac{1}{|W|\, \psi^2}{\sum}_r \xi^{b\, \frac{|\mu|^2
-|\rho|^2}{2}} \left( \sum_{w\in W} \sn(w) \xi^{(\mu|w(\rho))} \right)^2\notag \\
&=
\frac{1}{|W|\,  \psi^2}{\sum}_r \xi^{b\, \frac{|\mu|^2 -|\rho|^2}{2}} \sum_{w,w'\in W}
\sn(w w') \xi^{(\mu|w(\rho) + w'(\rho))} \notag \end{align}

Since $w(\rho) \in \rho+Y$ we have $w(\rho) + w'(\rho) \in 2\rho+Y =Y$. Using
(\ref{quadratic}),

\begin{align} F^{P\fg}_{U_b}(\xi)&= \frac{1}{|W| \psi^2} \gamma^{P\fg}_{b}(\xi)\,
\xi^{-b^*|\rho|^2} \sum_{w,w'\in W}\sn(w w')  \xi^{-b^* (w(\rho) | w'(\rho)) }\notag \\
&=
\frac{\gamma_b^{P\fg}\, \xi^{-b^* |\rho|^2}}{|W| \psi^2} |W| \sum_{w\in W} \sn(w)
\xi^{(-b^* \rho| w(\rho))}\notag\\
& = \frac{\gamma^{P\fg}_{b}(\xi)\,  \xi^{-b^*| \rho|^2}}{ \psi} J_U(-b^* \rho)\notag \\
&=  \frac{   \gamma_{b}^{P\fg}(\xi)\,     \xi^{(1-b^*) |\rho|^2}
J_U(b^*\rho)}{\prod_{\al>0} (1 -\xi^{(\al|\rho)})} \qquad \text{by (\ref{ssss}).}\notag
\end{align}

b)  follows from part a), Lemma \ref{Gammab}, and the fact that $s =(dim \fg -\ell)/2$.
\end{pf}

\begin{cor}  Suppose $r\ge d h^\vee$  is an odd prime, and $b$ is not divisible by $r$. Then
 $\tau^{P\fg}_M(\xi) \sim 1$, i.e. $\tau^{P\fg}_M(\xi)$ is a unit in $\Z[\xi]$,
  for the lens space $M$ obtained by surgery along $U_b$.
  \end{cor}

\begin{rem}
The actual value of $\tau_M^{P\fg}(\xi)$, where $M$ is obtained by surgery on $U_b$ is
(again here $b$ is an integer not divisible by the odd prime $r$)

\begin{equation}
  \tau_M^{P\fg}(\xi) = \left(\frac{|b|}{r} \right)^\ell \, \xi^{(\frac{\sn(b)
-b}{2}|\rho|^2)^{\sim}} \, \prod_{\al >0} \frac{1 - \xi^{-
(b^*\rho|\al)}}{1-\xi^{-(\sn(b)\rho|\al)}}.
\label{lens}
\end{equation}

Here  $\left(\frac{|b|}{r} \right)$ is the Legendre symbol, $(\frac{x}{y})^{\sim}$ is the
reduction modulo $r$, i.e.  $(\frac{x}{y})^{\sim} = x y^*$.
\end{rem}

\subsection{Expansion of quantum link invariants}

\begin{lem}\label{decomposition}
 For each positive integer $N$ one has

$$ Q_L(\mu_1,\dots,\mu_m) = \sum_{n=0}^{N-1}p_n(\mu_1,\dots,\mu_m) (q-1)^n + R,$$ where
$R$ is in $\Z[q^{\pm 1}]$ and divisible by $(q-1)^N$, $p_n(\mu_1,\dots,\mu_m)$ is a
polynomial function on $\hh^*_\R$ which takes integer values when $\mu_j\in Y$. Moreover
the degree of $p_n$ satisfies

\begin{equation}
deg (p_n) \le 2n + m(dim \fg -\ell) \label{degree}
\end{equation}
\end{lem}
\begin{pf}  This follows easily from a counting argument in
the theory of the Kontsevich integral, using the fact that $J_L$ is obtained from the
Kontsevich integral by substituting the Lie algebra into the chord diagrams (see
\cite{LM,Kassel}). The fact that $p_k$ takes integer values when $\mu_1,\dots,\mu_m\in Y$
follows from the integrality of the coefficients of $J_L$. Let us briefly sketch the 
idea. 

Expanding $J_L$ using $q= e^\hbar$, with $\hbar$ a new variable, we get

 $$J_L(\mu_1,\dots,\mu_m)|_{q= \exp \hbar} = \sum_{n=0}^\infty p'_{n}(\mu_1,\dots,\mu_m) \hbar^n,$$
 where $p'_n$ is a
function on $(\hh^*_\R)^m$. The Kontsevich integral theory will show that $p'_n$ is a
polynomial function with degree at most $2n + deg(dim(\mu_1)) +\dots+ deg( dim(\mu_m))$,
where $dim(\mu)$ is the function which gives the  dimension of the module
$\Lambda_{(\mu-\rho)}$. By Weyl formula, $dim(\mu)$ is a polynomial function of degree
$s$ -- the number of positive roots. Hence

$$ deg (p'_n) \le 2n + m s.$$

Thus for $Q_L= J_L J_{U^m}$ we have

$$ Q_L(\mu_1,\dots,\mu_m)|_{q= \exp \hbar} = \sum_{n=0}^\infty  p''_n(\mu_1,\dots,\mu_m)
\hbar^n,$$ where the degree of $p''_n$ is less than or equal to $2n + 2 ms$.

Change the variable from $\hbar$ to $q-1= e^\hbar-1$ (or $\hbar = \ln [(q-1) +1]$, we get

$$ Q_L(\mu_1,\dots,\mu_m) = \sum_{k=0}^\infty  p_n(\mu_1,\dots,\mu_m) (q-1)^n,$$ with
$deg (p_n) \le 2n + 2ms$. It remains to notice that s =($dim \fg - \ell)/2$.
\end{pf}

\subsection{A technical lemma}
\begin{lem} Suppose  $r$ is an odd prime, $p$ a polynomial
function taking values in $\Z$ when $\mu_1,\dots,\mu_m\in Y$. Let
$$ x= \sum_{\mu_j \in \rho + (P_r \cap Y)}  p(\mu_1,\dots,\mu_m)$$ and
 $$y = (\xi-1)^{\ell m
\frac{r-1}{2} - \lfloor \frac{deg p}{2}\rfloor },$$   where $\lfloor z\rfloor$ is the
greatest integer less than or equal to $z$. Then $x/y \in \Z[\xi]$. (Note that $y$ may
not be in $\Z[\xi]$.)
 \label{divisible}
\end{lem}

\begin{pf} In \cite[Corollary 4.14]{perturbative} we proved that the quotient $x/y$ is in $\Z[\xi,
\frac{1}{(r-1)!}]$. But  $y$ is coprime with $(r-1)!$, hence the quotient must be in
$\Z[\xi]$.
\end{pf}

\subsection{Proof of the theorem}

 Let $N= m \frac{r \ell -dim \fg}{2}$. Then the denominator of (\ref{definition3}) is a factor of
$(\xi-1)^N$ by Proposition \ref{ss23}. We will prove that the numerator is divisible by
$(\xi-1)^N$.

Applying Lemma \ref{decomposition}
$$ Q_{L}(\mu_1,\dots,\mu_m) =  \sum_{n=0}^{N-1}p_n(\mu_1,\dots,\mu_m) (q-1)^n + R,$$

We sum over $\mu_j\in \rho +(P_r\cap Y)$ to get $F_L^{P\fg}(\xi)$.  The term involving
$R$ is certainly divisible by $(\xi-1)^N$. For each $n$ the term involving $p_n$, by
Lemma \ref{divisible}  is divisible by

$$  (\xi-1)^{\ell m \frac{r-1}{2} - \lfloor \frac{\deg p_n}{2}\rfloor} \times (\xi-1)^n,$$ which, by
(\ref{degree}), is divisible by $(\xi-1)^N$.  This completes the proof of the theorem.

\section{Perturbative expansion}

\subsection{General }  Unlike the link case, quantum 3-manifold invariants can be defined only at
roots of unity, i.e. the domain of the function $\tau^\fg_M(q)$ is the set of rational
points on the unit circle in the complex plane $\Bbb C$. For many manifolds, eg the
Poincare sphere or the Brieskorn sphere $\Sigma(2,3,7)$, there is no analytic extension
of the function $\tau^\fg_M$ around $q=1$. In perturbative theory, we want to expand the
function $\tau^\fg_M$ around $q=1$ into power series. For rational homology 3-spheres,
i.e.manifolds $M$ with 0 rational homology, and for $\fg=sl_2$, Ohtsuki showed that there
is a number-theoretic expansion of $\tau^{P\fg}_M$ around $q=1$, see \cite{Ohtsuki}. We 
established the same result for the case $\fg  = sl_n$, see \cite{perturbative}. The 
proof in \cite{perturbative} is readily applied to any simple Lie algebra: In 
\cite{perturbative} we had to use some integrality properties of quantum link invariants 
and quantum 3-manifold invariants, and there  we established these properties for the 
special case $\fg= sl_n$. For the general simple Lie algebras, these integrality 
properties are the results of \cite{integrality} and Theorem \ref{integ3}.

\subsection{The number-theoretic expansion}
 Suppose $r$ is a big enough prime, and
$\xi = \exp (2 \pi i/r)$. By the integrality (Theorem \ref{integ3}),

$$\tau^{P\fg}_M(\xi) \in \Z[\xi] =  \Z[q] /(1 +q + q^2 + \dots + q^{r-1}).$$

Choose a representative $f(q)\in \Z[q]$ of $\tau^{P\fg}_M(\xi)$. Formally substitute
$q=e^\hbar$ in $f(q)$:

$$f(q) = c_{r,0} + c_{r,1} \hbar + \dots+  c_{r,n} \hbar^n + \dots$$

The rational numbers $c_{r,n}$ depend on $r$ and the representative $f(q)$. Their
denominators must be a factors of $n!$, by Theorem \ref{integ3}. Hence if $n<r-1$, we can 
reduce $c_{r,n}$ modulo $r$ and get an element of $\Z/r\Z$. It is easy to see that these 
reductions $c_{r,n} \pmod r$ do not depend on the representative $f(q)$ and hence are 
invariants of the 3-manifolds. The dependence on $r$ is a big drawback. The theorem below 
says that there is a number $c_n$, not depending on $r$, such that $c_{r,n} \pmod r$ is 
the reduction of $c_n$, or $-c_n$, modulo $r$, for sufficiently large prime $r$. It is 
easy to see that if such $c_n$ exists, it must be unique. 

\begin{thm} For every ratiomal homology 3-sphere $M$, there are a sequence of numbers $c_n \in \Z[\frac{1}{(2n+ 2s)!
|H_1(M,\Z)|}]$, such that for sufficiently large prime $r$ (actually any prime $r> \max
(|H_1(M,\Z), dim \fg -\ell)$ is enough),

$$ c_{r,n}  \equiv \left( \frac{|H_1(M,\Z)|}{r}\right)^\ell \, c_n \pmod r,$$ where
$\left(\frac {|H_1(M,\Z)|}{r}\right)=\pm 1 $ is the Legendre symbol.
\end{thm}

The series  ${\frak t}^{P\fg}_M(\hbar) = \sum_{n=0}^\infty c_n \hbar^n$ can be considered
as the perturbative expansion of the function $\tau_M^{P\fg}$ at $q=1$.  As mentioned
above, the proof is just similar to the one for the case $\fg =sl_n$ in
\cite{perturbative}.

\subsection{Some calculation}  Let us describe here how to calculate the power series ${\frak
t}_M$, and sketch the ideas behind the calculation.

\subsubsection{The $\fg=sl_2$  and surgery on a knot case}  In this case let the positive integer $N$ stand for the
unique $\fg$-module of dimension $N$. The invariant $J_L(N_1,\dots,N_m)$ is known as the
colored Jones polynomial. Suppose $M$ is obtained by surgery along a knot $K$ with
framing 1. Let $K^0$ be the same knot with framing 0.  Then

\begin{equation} Q_{K^0}(N)|_{q = e^{\hbar}} = \sum_{2 \le j \le n+2} c_{j,n} N^{j} \hbar^n.
\label{ss}
\end{equation}
 The restriction $2j\le n+2$ follows from the fact that $K^0$ has framing 0. It is known
 that there is no odd order of $N$: $j$ must be even.

To obtain ${\frak t}_M(\hbar)$, all one needs is to replace $N^{2j}$ in (\ref{ss}) by
$(-2)^j(2j-1)!!\, \hbar^{-j}$, then multiply by a universal constant:

$$ {\frak t}_M(\hbar) = z \sum c_{2j,n} (-2)^j (2j-1)!!\, \hbar^{n-j},$$
where $z= z^{sl_2} = (1-q)/2 = (1-e^\hbar)/2$.

Presumably this was first obtained by Rozansky \cite{Rozansky}.

\subsubsection{The case of general simple Lie algebra and surgery along a knot}

Again assume that $M$ is obtained by surgery along the knot $K$ with framing 1, and $K^0$
is the same knot with framing 0. It is known that every polynomial function $p(\mu)$ on
$\hh^*_\R$ are linear combinations of functions of the form $\beta^j$, $\beta \in Y$.
Here $\beta^j(\mu) := (\beta|\mu)^j$. Thus one has

\begin{equation}
Q_{K^0}(\mu)|_{q = e^\hbar} = \sum_{2s \le j \le n+2s,\, \beta \in Y} c_{\beta;j;n}
\beta^j(\mu)\, \hbar^n.
\label{ss1}
\end{equation}

Here for each degree $n$ the sum is finite. Again the restriction $j \le 2s +n$ comes
from the fact that $K^0$ has framing 0.

To obtain ${\frak t}_M(h)$, all one needs is to replace $\beta^j(\mu)$ in (\ref{ss1}) by
0 if $j$ is odd, $\beta^{2j}(\mu)$ by
\begin{equation}(2j-1)!!\,\hbar^{-j}\, (-|\beta|^2)^j, \label{tt}
\end{equation}
then multiply by a universal constant:

\begin{equation} {\frak t}^{P\fg}_M(\hbar) =  \sum c_{\beta;2j;n} (2j-1)!!\, (-|\beta|^2)^j\, \hbar^{n-j} \, \times
\, \frac{1}{|W|} \prod_{\al >0} (1-q^{(\al|\rho)})
\end{equation}

\subsubsection{A sketch of the main idea} The main idea is to separate the framing part, and consider the sum
${\sum}_r$ as a discreet Gauss integral. This was first used by Rozansky (for $sl_2$) in
his series of important work on quantum invariants.

 Recall how we define $\tau^{P\fg}_M(\xi)$.
One gets $Q_K$ by multiplying $Q_{K^0}$ by $q^{(|\mu|^2 -|\rho|^2)/2}$. Summing $Q_K$
over $\mu\in \rho + (P_r\cap Y)$, we get $F_K$. Then we have to divide $F_K$ by 
$F_{U_+}$. The result is $\tau^{P\fg}_M(\xi)$. A look at formula (\ref{ss1}) shows that 
if we understand the perturbative expansion of 
 \begin{equation} \frac{{\sum}_r
q^{\frac{|\mu|^2-|\rho|^2}{2}} \beta^j(\mu)}{F_{U_+}},\label{uu}
\end{equation} then we will know
the perturbative expansion  of $\tau_M^{P\fg}$.

If we replace $\beta^j(\mu)=(\beta|\mu)^j$ in (\ref{uu}) by $q^\beta(\mu):=
q^{(\beta|\mu)}$, then the perturbative expansion  is easy to calculate:

\begin{align}
 \frac{{\sum}_r \xi^{\frac{|\mu|^2-|\rho|^2}{2}} \xi^{(\beta|\mu)}}{F_{U_+}}
 &= \frac{1}{F_{U_+}} \gamma^{P\fg}(\xi) \, \xi^{-|\beta|^2/2} \qquad \text{by (\ref{quadratic})} \label{ss3} \\
 &= \xi^{-|\beta|^2/2}\, \prod_{\al>0}(1-\xi^{(\al|\rho)}) \qquad \text{by (\ref{x2})}\notag \end{align}
Thus the perturbative expansion  of the left hand side of (\ref{ss3}) is
$q^{-|\beta|^2/2} z$, with $z= \prod_{\al>0} (1-q^{(\al|\rho)})$ and $q= e^\hbar$.

Now if we expand $q^{(\beta|\mu)}= \exp[\hbar(\beta|\mu)]$, we can see the term
$(\beta|\mu)^j$ there:

$$ \exp[\hbar (\beta|\mu)] = \sum_{j \ge 0} \frac{\hbar^j (\beta|\mu)^j}{j!}.$$

To obtain the perturbative expansion of (\ref{uu}), we expand $q^{-|\beta|^2/2} z$ into
power series of $\hbar$ , and keep only the part of degree $j$ in $\mu$. It is easy to
see that if $j$ is odd, there is no part of degree $j$, and if $j$ is even, then the part
of degree $j$ is given by the formula (\ref{tt}) (In this argument we consider $\mu$ as a
variable. To be more precise, one replace $\mu$ by $t \mu$, with $t\in \R$ a variable,
then compare the terms of same degree of $t$.)

\subsubsection{Special lens spaces} Suppose $M$ is obtained by surgery on the unknot
$U_b$, with $b\neq 0$. Then  from  (\ref{lens}) it follows that

$${\frak t}^{P\fg}_M(\hbar) =  q^{\frac{\sn(b) -b}{2}|\rho|^2} \prod_{\al>0} \frac{1- q^{-(\rho|\al)/b}}{1-
q^{-\sn(b)(\rho|\al)}}\mid_{q =\exp \hbar}.$$

\subsubsection{Link with diagonal linking matrix} Suppose $L$ is a framing link whose linking matrix is
diagonal, with non-zero integers $b_1,\dots,b_m$ on the diagonal. Let $L^0$ be the same
link with 0 framing, and $M$ the 3-manifold obtained by surgery along $L$, which is a
rational homology 3-sphere. Expanding $q=e^\hbar$ in $Q_{L^0}$ we get

$$ Q_{L^0}(\mu_1,\dots,\mu_m)|_{q= e^\hbar} = \sum _{\beta_1,\dots,\beta_m\in Y;\,j_1,\dots,j_m\in Z_+;\,n\in \Z_+}
c_{\beta_1,\dots,\beta_m;j_1,\dots,j_m;n}\,\beta_1^{j_1}(\mu_1) \dots 
\beta_m^{j_m}(\mu_m) \, \hbar^n.$$ 

There are some restrictions on $\beta_i,j_i$, for a fixed $n$. Then to obtain ${\frak
t}^{P\fg}_M(\hbar)$ one needs to replace $\beta_i^{j}(\mu_i)$ by 0 if $j$ is odd,
$\beta_i^{2 j}(\mu_i)$ by
$$z_{b_i}\, b_i^{-j}\,  (2j-1)!!\, (-|\beta_i|^2)^j\, \hbar^{-j},$$
where
$$z_{b_i} = \frac{1}{|W|} \, q^{\frac{|\rho|^2}{2}(\sn(b) -b)}\prod_{\al >0} (1-q^{\sn (b)
(\al|\rho)}).$$

Thus,

$${\frak t}^{P\fg}_M(\hbar) = z_{b_1} \dots z_{b_m} \sum c_{\beta_1,\dots,\beta_m;2j_1,\dots,2j_m;n}
\prod_{i=1}^m  (2j_i-1)!!\, (\frac{-|\beta_i|^2}{b_i})^{j_i}\, \hbar^{n - j_1-\dots
-j_m}.$$

The restriction on $j_1,\dots,j_m$ will guarantee that the right hand side is a formal
power series in $h$.

\subsubsection{General case}  Suppose now $M$ is an arbitrary rational homology 3-sphere.
Ohtsuki showed that there are lens spaces $M_1,\dots,M_l$, each obtained by surgery on an
unknot with non-zero framing, such that $M'=M\# M_1 \# \dots \# M_l$ can be obtained
surgery along a link with diagonal linking matrix, see \cite{Ohtsuki}. Then one has

$$ {\frak t}^{P\fg}_M(\hbar) =  {\frak t}^{P\fg}_{M'}(\hbar) \, \left( {\frak
t}^{P\fg}_{M_1}(\hbar)\right)^{-1} \dots  \left( {\frak
t}^{P\fg}_{M_1}(\hbar)\right)^{-1}.$$

        \end{document}